# AN APPLICATION OF LIMITING INTERPOLATION TO THE FOURIER SERIES THEORY

## Leo R. Ya. Doktorski

Department Object Recognition, Fraunhofer Institute of Optronics, System Technologies and Image Exploitation IOSB, Gutleuthausstr. 1, 76275 Ettlingen, Germany

**Abstract.** Limiting real interpolation method is applied to describe the behaviour of the Fourier coefficients of functions that belong to spaces which are "very close" to $L_2$.

## 1.   Introduction

We consider (equivalent classes of) complex-valued measurable functions on [0,1]. Let $\{\varphi_n\}$ ($n \in \mathbb{N}$ or $n \in \mathbb{Z}$) be an orthonormal system in $L_2$ bounded in $L_\infty$:

$$\sup \| \varphi_n \mid L_\infty \| = M \; (<\infty). \tag{1.1}$$

Everywhere below we denote by $\{c_n(f)\}$ the Fourier coefficients of a function $f$ with respect to the system $\{\varphi_n\}$:

$$c_n \equiv c_n(f) := \int_0^1 f(x)\overline{\varphi_n}(x)\,dx.$$

We write $\mathcal{F}$ for the Fourier series map assigning the sequence of Fourier coefficients to any function $f$, i.e. $\mathcal{F}(f) = \{c_n(f)\}$. It is known that $\mathcal{F}$ is a linear bounded operator from $L_2$ to $l_2$ and from $L_1$ to $l_\infty$:

$$\| \mathcal{F} \mid L_2 \to l_2 \| = 1, \tag{1.2}$$

$$\| \mathcal{F} \mid L_1 \to l_\infty \| \leq M. \tag{1.3}$$

The results of Hausdorff, Young and Paley describe the behaviour of the Fourier coefficients of functions that belong to Lebesgue spaces $L_p$ (1<p<2). For further information about classical results dealing with Fourier series map we refer to e.g. [2, 3, 27, 28, 39]. Interpolation between (1.2) and (1.3) by the classical real interpolation method $(*,*)_{\theta,q}$ (0<$\theta$<1) provides such description for the Lorentz spaces $L_{p,q}$ [1, 35]:

$$\| \mathcal{F} \mid L_{p,q} \to l_{r,q} \| \prec M^{\frac{2}{p}-1}.$$

Application of the real interpolation functor $(*,*)_{\theta,q,\alpha}$ (0<$\theta$<1, $\alpha \in \mathbb{R}$) involving logarithmic factor shows that for the Lorentz–Zygmund spaces $L_{p,q}(logL)_\alpha$ it holds [2, 14, 15]

$$\mathcal{F} \colon L_{p,q}(logL)_\alpha \to l_{r,q}(\log l)_\alpha.$$

*E-mail address*: leo.doktorski@iosb.fraunhofer.de, *Phone number*: +49 7243 992-217.

In both formulae $1<p<2$, $\frac{1}{r}=1-\frac{1}{p}$, $0<q\leq\infty$. Note that in the scale $L_{p,q}(logL)_\alpha$ this result is optimal [34, Theorem 5.3]. Similar results were obtained also by other approaches (see [26, 33, 34] and references therein).

But for the spaces which are "very close" to $L_2$ all these methods do not work and two other approaches can be applied. One of them is a "direct way" based on estimates of a K-functional [5, 6, 29, 30, 32]. The other one is based on limiting real interpolation methods (see e. g. [7, 14, 16, 36]).

Comparable results are also known for the inverse transformation

$$\mathcal{F}^{-1}: \{c_n\} \to f(x):=\sum_{n=1}^{\infty} c_n \varphi_n(x).$$

Here $\{c_n\}$ are sequences of complex numbers. It is known [35, 39] that

$$\|\mathcal{F}^{-1} \mid l_2 \to L_2\| = 1, \tag{1.4}$$

$$\|\mathcal{F}^{-1} \mid l_1 \to L_\infty\| = M. \tag{1.5}$$

Interpolation between (1.4) and (1.5) by the real interpolation method $(*,*)_{\theta,q}$ or $(*,*)_{\theta,q,\alpha}$ ($0<\theta<1$) leads also to results for Lorentz and Lorentz–Zygmund sequence spaces with main parameter laying between 1 and 2 [14, 15, 35]. But for the spaces which are "very close" to $l_2$ all these methods do not work either.

The main objective of this work is to study the Fourier series map and its inverse for the spaces which are "very close" to $L_2$ or $l_2$ resp. via limiting interpolation methods. More precisely, we prove the following assertions.

**Theorem 1.1.** Let $0< q \leq \infty$ and $\alpha< -1/q$. Then

$$\left(\sum_{k=1}^{\infty}\left[(1+|\log k|)^\alpha \left(\sum_{i=1}^{k}\left(c_i^*(f)\right)^2\right)^{1/2}\right]^q k^{-1}\right)^{1/q} \leq$$

$$\leq C \min(M, (1+\log M)^{|\alpha|}) \left(\int_0^1\left[(1+|\log t|)^\alpha \left(\int_t^1 \left(f^*(u)\right)^2 du\right)^{1/2}\right]^q t^{-1} dt\right)^{1/q}$$

for some constant $C$ which depends only on $q$ and $\alpha$. (As usual, the integral and the sum should be replaced by the supremum when $q=\infty$.)

**Theorem 1.2**. Let $0<q\leq \infty$ and $\alpha> -1/q$. Then for any $\varepsilon>0$

$$\left(\int_0^1\left[(1+|\log t|)^\alpha \left(\int_0^t \left(f^*(u)\right)^2 du\right)^{1/2}\right]^q t^{-1} dt\right)^{1/q} \leq$$

$$\leq C \min(M, (1+\log M)^{|\alpha|+\varepsilon+1/q}) \left( \sum_{k=1}^{\infty} \left[ (1+|\log k|)^{\alpha} \left( \sum_{i=k}^{\infty} \left(c_i^*(f)\right)^2 \right)^{1/2} \right]^q k^{-1} \right)^{1/q}$$

for some constant $C$ which depends only on $q$, $\alpha$, and $\varepsilon$. (As usual, the integral and the sum should be replaced by the supremum when $q=\infty$.) In particular, this means that if the expression on the right-hand side of the estimate above exists then the row $\sum_{n=1}^{\infty} c_n \varphi_n(x)$ converges in the (quasi-) norm determined by the expression on the left-hand side.

**Remark 1.2.**

(i) The orthonormality of the system $\{\varphi_n\}$ and (1.1) implies $M \geq 1$.

(ii) The system $\{\varphi_n\}$ may also be bounded in $L_Q$ for some $Q \in (2,\infty]$ [28]. The estimates of Theorems 1.1 and 1.2 do not depend on $Q$.

The second objective of the present paper is to compare the estimate given in Theorem 1.1 with known results. It turns out that the results established with the help of the "direct way" approach are weaker than those of Theorem 1.1.

This paper is organized as follows. Section 2 contains necessary notations, definitions and auxiliary results. Theorems 1.1 and 1.2 will be proven and discussed in Sections 3 and 4 resp.

## 2. Notation, definitions and auxiliary results

If $X$ is a (quasi-) Banach space and $x \in X$ then its (quasi-) norm is denoted as $\| x \mid X \|$. $X \cong Y$ means that the Banach spaces $X$ and $Y$ are isomorphic. Throughout the paper $L_q(a,b)$ ($0 < q \leq \infty$, $-\infty \leq a < b \leq \infty$) is the usual quasinormed Lebesgue space $L_q$ on the interval $(a,b)$. For $q \geq 1$ it is a Banach space. $L_q$ implies $L_q(0,1)$. By $C$ we designate different positive constants wich are independent of all significant arguments. If $f$ and $g$ are positive functions, we will write $f \prec g$ if $f \leq C \cdot g$ and $f \approx g$ if $f \prec g$ and $g \prec f$.

### 2.1. *Interpolation spaces*

Let $\vec{X} = (X_0, X_1)$ be a compatible couple of (quasi-) Banach spaces and let

$$K(t,x) \equiv K(t,x,\vec{X}) := \inf_{x = x_0 + x_1; x_i \in X_i} (\|x_0 \mid X_0\| + t\|x_1 \mid X_1\|).$$

be Peetre's $K$-functional. For further information about properties of the $K$-functional and the real interpolation method, we refer to [3, 4, 27, 38]. For our purposes, it is enough to consider only ordered couples $X_0 \supset X_1$ with the norm of embedding equal to 1. This will be denoted as $X_0 \overset{1}{\supset} X_1$. In this case $K(t,x) \approx \|x \mid X_0\|$ for $t > 1$ [4].

**Definition 2.1.** Let $X_0 \overset{1}{\supset} X_1$, $0 \leq \theta \leq 1$, $0 < q \leq \infty$, and $\alpha \in \mathbb{R}$. We set

$$\vec{X}_{\theta,q,\alpha} \equiv (X_0, X_1)_{\theta,q,\alpha} := \{x \in X_0 + X_1 \mid \|x \mid \vec{X}_{\theta,q,\alpha}\| := \|t^{-\theta-1/q}(1+|\log t|)^\alpha K(t,x) \mid L_q(0,1)\| < \infty\}.$$

It only makes sense to consider the spaces $\vec{X}_{\theta,q,\alpha}$ on the set

$$\{(\theta,q,\alpha) \in [0, 1] \times (0, \infty] \times \mathbb{R} \mid 0 < \theta < 1, \text{ or } \theta = 0, \, q \leq \infty, \, \alpha \geq -1/q,$$
$$\text{or } \theta = 1, \, q \leq \infty, \, \alpha < -1/q, \text{ or } \theta = 1, \, q = \infty, \, \alpha = 0\}.$$

Note that the functors $(X_0, X_1)_{0,q,\alpha}$ and $(X_0, X_1)_{1,q,\alpha}$ produce spaces which are "very close" to $X_0$ and to $X_1$ respectively. This definition can be found in a lot of papers (see e. g. [10, 12, 15, 16, 17, 19, 21, 36]).

Different analogues of the next lemma can be found in literature. See [12, Theorem 2.5], [8, Theorem 4.9], [9, (2.3) and (2.4)], [18], and [22, Theorem 3.5].

**Lemma 2.2.** Let $X_0 \overset{1}{\supset} X_1$ and $Y_0 \overset{1}{\supset} Y_1$ be (quasi-) Banach spaces, and let $T$ be a (quasi-) linear bounded operator, $T: X_j \to Y_j$ with the norms $M_j := \|T \mid X_j \to Y_j\|$ ($j = 0, 1$). Additionally suppose that $0 < q \leq \infty$.

(a) If $M_0 \geq M_1$ and $\alpha < -1/q$, then $T$ is bounded from $\vec{X}_{1,q,\alpha}$ to $\vec{Y}_{1,q,\alpha}$ and

$$\|T \mid \vec{X}_{1,q,\alpha} \to \vec{Y}_{1,q,\alpha}\| \leq \min(M_0, (1+\log(M_0/M_1))^{|\alpha|} M_1).$$

(b) If $M_0 \leq M_1$ and $\alpha \geq -1/q$, then $T$ is bounded from $\vec{X}_{0,q,\alpha}$ to $\vec{Y}_{0,q,\alpha}$ and for any $\varepsilon > 0$

$$\|T \mid \vec{X}_{0,q,\alpha} \to \vec{Y}_{0,q,\alpha}\| \prec \min(M_1, (1+\log(M_1/M_0))^{|\alpha|+\varepsilon+1/q} M_0).$$

**Proof.** First, notice that if $x \in X_0 + X_1$ then

$$K(t, Tx; \vec{Y}) \leq \max(M_0, M_1) K(t, x; \vec{X}). \tag{2.1}$$

Moreover,

$$K(t, Tx; \vec{Y}) \leq M_0 K(tM_1/M_0, x; \vec{X}). \tag{2.2}$$

We begin with the assertion (a). (2.1) implies

$$\|T \mid \vec{X}_{1,q,\alpha} \to \vec{Y}_{1,q,\alpha}\| \leq M_0. \tag{2.3}$$

It is not difficult to show that

$$(1+|\log(uv)|)^\alpha \leq (1+|\log u|)^\alpha (1+|\log v|)^{|\alpha|} \qquad (u,v > 0).$$

By means of this inequality and of (2.2), because $M_1/M_0 \leq 1$, and using the change of variable $u = tM_1/M_0$, we obtain

$$\| Tx \mid \vec{Y}_{1,q,\alpha} \| = \left\| t^{-1-1/q}(1+|\log t|)^\alpha K(t,Tx;\vec{Y}) \mid L_q(0,1)\right\| \le$$

$$\le M_0 \left\| t^{-1-1/q}(1+|\log t|)^\alpha K(tM_1/M_0, x;\vec{X}) \mid L_q(0,1)\right\| =$$

$$= M_1 \left\| u^{-1-1/q}(1+|\log(uM_0/M_1)|)^\alpha K(u,x;\vec{X}) \mid L_q(0, M_1/M_0)\right\| \le$$

$$\le M_1(1+\log(M_0/M_1))^{|\alpha|} \left\| u^{-1-1/q}(1+|\log u|)^\alpha K(u,x;\vec{X}) \mid L_q(0,1)\right\| =$$

$$= M_1(1+\log(M_0/M_1))^{|\alpha|} \| x \mid \vec{X}_{1,q,\alpha} \|.$$

Combining this with (2.3) we get the assertion (a). Next we consider the case (b). (2.1) implies

$$\| T \mid \vec{X}_{0,q,\alpha} \to \vec{Y}_{0,q,\alpha} \| \le M_1. \tag{2.4}$$

For real numbers $\alpha_0$ und $\alpha_\infty$ we use as usual

$$l^{(\alpha_0,\alpha_\infty)}(t) := \begin{cases} (1-\log t)^{\alpha_0}, & \text{if } 0 < t \le 1, \\ (1+\log t)^{\alpha_\infty}, & \text{if } 1 \le t < \infty. \end{cases}$$

Now we show that for any $\varepsilon > 0$

$$\left\| t^{-1/q}(1+|\log t|)^\alpha K(t,Tx;\vec{Y}) \mid L_q(0,1)\right\| \approx \left\| t^{-1/q} l(t)^{(\alpha,-\varepsilon-\frac{1}{q})} K(t,Tx;\vec{Y}) \mid L_q(0,\infty)\right\|. \tag{2.5}$$

It is enough to check that

$$\left\| t^{-1/q}(1+|\log t|)^{-(\varepsilon+\frac{1}{q})} K(t,Tx;\vec{Y}) \mid L_q(1,\infty)\right\| \prec \left\| t^{-1/q}(1+|\log t|)^\alpha K(t,Tx;\vec{Y}) \mid L_q(0,1)\right\|.$$

Observe that

$$\left\| t^{1-1/q}(1+|\log t|)^\alpha \mid L_q(0,1)\right\| < \infty$$

and if $\varepsilon > 0$

$$\left\| t^{-1/q}(1+|\log t|)^{-(\varepsilon+\frac{1}{q})} \mid L_q(1,\infty)\right\| < \infty.$$

Using that $t^{-1}K(t,x)$ is non-increasing and $K(t,Tx;\vec{Y}) \approx K(1,Tx;\vec{Y})$ for $t>1$, we obtain

$$\left\| t^{-1/q}(1+|\log t|)^{-(\varepsilon+\frac{1}{q})} K(t,Tx;\vec{Y}) \mid L_q(1,\infty)\right\| \approx K(1,Tx;\vec{Y}) \left\| t^{-1/q}(1+|\log t|)^{-(\varepsilon+\frac{1}{q})} \mid L_q(1,\infty)\right\| \approx$$

$$\approx K(1,Tx;\vec{Y}) \approx \frac{K(1,Tx;\vec{Y})}{1} \left\| t^{1-1/q}(1+|\log t|)^\alpha \mid L_q(0,1)\right\| \le \left\| t^{-1/q}(1+|\log t|)^\alpha K(t,Tx;\vec{Y}) \mid L_q(0,1)\right\|.$$

So, (2.5) is proven. It can be shown that (cf. [9], p. 169.)

$$l^{(\alpha_0,\alpha_\infty)}(uv) \le l^{(\alpha_0,\alpha_\infty)}(u)(1+|\log v|)^{|\alpha_0|+|\alpha_\infty|} \qquad (u,v>0).$$

Therefore,

$$l^{(\alpha,\,-\varepsilon-1/q)}(uv) \le l^{(\alpha,\,-\varepsilon-1/q)}(u)(1+|\log v|)^{|\alpha|+\varepsilon+1/q}.$$

Since $M_0 \leq M_1$, using (2.2) and (2.5), and by means of the change of variable $u=tM_1/M_0$, we obtain:

$$\| Tx \mid \vec{Y}_{0,q,\alpha} \| = \left\| t^{-1/q}(1+|\ln t|)^\alpha K(t,Tx;\vec{Y}) \mid L_q(0,1) \right\| \approx$$

$$\approx \left\| t^{-1/q} l(t)^{(\alpha,-\varepsilon-\frac{1}{q})} K(t,Tx;\vec{Y}) \mid L_q(0,\infty) \right\| \leq$$

$$\leq M_0 \left\| t^{-1/q} l(t)^{(\alpha,-\varepsilon-\frac{1}{q})} K(tM_1/M_0,x;\vec{X}) \mid L_q(0,\infty) \right\| =$$

$$= M_0 \left\| u^{-1/q} l(uM_0/M_1)^{(\alpha,-\varepsilon-\frac{1}{q})} K(u,x;\vec{X}) \mid L_q(0,\infty) \right\| \prec$$

$$\prec M_0(1+\log(M_1/M_0))^{|\alpha|+\varepsilon+1/q} \left\| u^{-1/q} l(u)^{(\alpha,-\varepsilon-\frac{1}{q})} K(u,x;\vec{X}) \mid L_q(0,\infty) \right\| \approx$$

$$\approx M_0(1+\log(M_1/M_0))^{|\alpha|+\varepsilon+1/q} \left\| u^{-1/q}(1+|\log u|)^\alpha K(u,x;\vec{X}) \mid L_q(0,1) \right\| =$$

$$= M_0(1+\log(M_1/M_0))^{|\alpha|+\varepsilon+1/q} \| x \mid \vec{X}_{0,q,\alpha} \|.$$

Combining this with (2.4) we complete the proof. $\square$

**Remark 2.3.** Due to Lemma 2.2 we have in Theorems 1.1 and 1.2 expression of the form
$$G(M,\gamma):=\min(M, (1+\log M)^\gamma)$$
with $M\geq 1$ and $\gamma\geq 0$. It is clear that $G(1,\gamma)=1$ and
$$G(M,\gamma) = \begin{cases} (1+\log M)^\gamma, & \text{if } \gamma \leq \log M / \log(1+\log M), \\ M, & \text{otherwise.} \end{cases}$$
In particular, $G(M,\gamma) = (1+\log M)^\gamma$ if $\gamma \leq 1$.

Next we consider the following limiting interpolation spaces. These spaces allow to formulate reiteration theorems in the limiting cases $\theta = 0$ and $\theta = 1$. In more general form these spaces were introduced and investigated in [14, 16, 19, 23, 25].

**Definition 2.4.** Let $X_0 \stackrel{1}{\supset} X_1$, $0<q,r \leq \infty$, $\alpha \in \mathbb{R}$. We denote by $\vec{X}^{\mathcal{L}}_{\theta,q,\alpha,r}$ ($0\leq\theta<1$) and $\vec{X}^{\mathcal{R}}_{\theta,q,\alpha,r}$ ($0<\theta\leq 1$) the sets of elements $x\in X_0$ for which the expessions
$$\left\| x \mid \vec{X}^{\mathcal{L}}_{\theta,q,\alpha,r} \right\| := \left\| t^{-1/q}(1+|\ln t|)^\alpha \left\| u^{-\theta-1/r} K(t,x) \mid L_r(0,t) \right\| \mid L_q(0,1) \right\|,$$

$$\left\| x \mid \vec{X}^{\mathcal{R}}_{\theta,q,\alpha,r} \right\| := \left\| t^{-1/q}(1+|\ln t|)^\alpha \left\| u^{-\theta-1/r} K(t,x) \mid L_r(t,1) \right\| \mid L_q(0,1) \right\|$$

resp. are finite.

The next lemma follows from [14, Lemma 6.2], [16, Lemma 4], see also [19].

**Lemma 2.5.** Let $0 < q,r \leq \infty$, $0\leq\theta<1$ and $\alpha>-1/q$ (or $0<\theta\leq 1$ and $\alpha<-1/q$); then

$$\vec{X}_{\theta,q,\alpha+\frac{1}{\min(r,q)}} \subset \vec{X}^{\mathcal{L}}_{\theta,q,\alpha,r} \text{ (or } \vec{X}^{\mathcal{R}}_{\theta,q,\alpha,r} \text{ resp.)} \subset \vec{X}_{\theta,r,\alpha+\frac{1}{\max(r,q)}} \cap \vec{X}_{\theta,\max(r,q),\alpha+\frac{1}{q}}.$$

## *2.2. Function spaces*

In this section we give necessary definitions of function and sequence spaces. We consider (equivalent classes of) complex-valued measurable functions on [0,1] and bounded complex-valued sequences $\{c_k\}$. As usual, $f^*$ is the non-increasing rearrangement of $|f|$, and $\{c_k^*\}$ ($k \in \mathbb{N}$) is the non-increasing rearrangement of the sequence $\{|c_k|\}$. Lorentz–Zygmund spaces can be defined as follows.

**Definition 2.6.** Let $0 < p,q \leq \infty$ and $\alpha \in \mathbb{R}$. Put

$$L_{p,q}(\log L)_\alpha := \{ f \mid \| f \mid L_{p,q}(\log L)_\alpha \| := \| t^{1/p-1/q}(1+|\log t|)^\alpha f^*(t) \mid L_q(0,1) \| < \infty \}.$$

Analogously

$$l_{p,q}(\log l)_\alpha := \{ \{c_k\} \mid \| \{c_k\} \mid l_{p,q}(\log l)_\alpha \| := \| k^{1/p-1/q}(1+\log k)^\alpha c_k^* \mid l_q \| < \infty \}.$$

These spaces are studied in [2, 13, 15, 16, 17, 19, 20, 21, 31]. See also [3] and [25]. Note that $L_{p,q} = L_{p,q}(\log L)_0$ and $l_{p,q} = l_{p,q}(\log l)_0$. Concerning next definition we refere to [14, 19, 25].

**Definition 2.7.** Let $0 < p,q,r \leq \infty$ and $\alpha \in \mathbb{R}$. Put

$$L^{\mathcal{L}}_{p,q,\alpha,r} := \{ f \mid \| f \mid L^{\mathcal{L}}_{p,q,\alpha,r} \| := \| t^{-1/q}(1+|\ln t|)^\alpha \| u^{1/p-1/r} f^*(u) \mid L_r(0,t) \| \mid L_q(0,1) \| \},$$

$$L^{\mathcal{R}}_{p,q,\alpha,r} := \{ f \mid \| f \mid L^{\mathcal{R}}_{p,q,\alpha,r} \| := \| t^{-1/q}(1+|\ln t|)^\alpha \| u^{1/p-1/r} f^*(u) \mid L_r(t,1) \| \mid L_q(0,1) \| \}.$$

Analogously

$$l^{\mathcal{L}}_{p,q,\alpha,r} := \left\{ \{c_k\} \mid \| \{c_k\} \mid l^{\mathcal{L}}_{p,q,\alpha,r} \| := \left( \sum_{k=1}^{\infty} \left[ k^{-\frac{1}{q}}(1+\ln k)^\alpha \left( \sum_{i=1}^{k} \left[ i^{\frac{1}{p}-\frac{1}{r}} c_i^* \right]^r \right)^{1/r} \right]^q \right)^{1/q} \right\},$$

$$l^{\mathcal{R}}_{p,q,\alpha,r} := \left\{ \{c_k\} \mid \| \{c_k\} \mid l^{\mathcal{R}}_{p,q,\alpha,r} \| := \left( \sum_{k=1}^{\infty} \left[ k^{-\frac{1}{q}}(1+\ln k)^\alpha \left( \sum_{i=k}^{\infty} \left[ i^{\frac{1}{p}-\frac{1}{r}} c_i^* \right]^r \right)^{1/r} \right]^q \right)^{1/q} \right\}.$$

(As usual, the sum should be replaced by the supremum when $q = \infty$ or $r = \infty$.)

We need only the combination $p = r = 2$ for these spaces:

$$\| f \mid L^{\mathcal{L}}_{2,q,\alpha,2} \| = \| t^{-1/q}(1+|\ln t|)^\alpha \| f^*(u) \mid L_2(0,t) \| \mid L_q(0,1) \|,$$

$$\|f \mid L^{\mathcal{R}}_{2,q,\alpha,2}\| = \|t^{-1/q}(1+|\ln t|)^{\alpha}\|f^{*}(u) \mid L_2(t,1)\| \mid L_q(0,1)\|,$$

$$\|\{c_k\} \mid l^{\mathcal{L}}_{2,q,\alpha,2}\| = \left(\sum_{k=1}^{\infty}\left[k^{-\frac{1}{q}}(1+\ln k)^{\alpha}\left(\sum_{i=1}^{k}\left(c_i^*\right)^2\right)^{1/2}\right]^q\right)^{1/q},$$

$$\|\{c_k\} \mid l^{\mathcal{R}}_{2,q,\alpha,2}\| = \left(\sum_{k=1}^{\infty}\left[k^{-\frac{1}{q}}(1+\ln k)^{\alpha}\left(\sum_{i=k}^{\infty}\left(c_i^*\right)^2\right)^{1/2}\right]^q\right)^{1/q}.$$

Note that in the terminology of [24] $L^{\mathcal{R}}_{2,\infty,\alpha,2}$ is the generalized grand Lorentz space $L^{2),2}_{(1+|\ln t|)^{\alpha}}$ and $l^{\mathcal{L}}_{2,\infty,\alpha,2}$ is the generalized grand Lorentz space of sequences $l^{2),2}_{(1+|\ln t|)^{\alpha}}$.

The following result is a consequence of [14, Corollaries 7.3 and 7.9], [19, Theorem 8.9], and [25, Theorem 5.7].

**Lemma 2.8.** Let $0<q\leq\infty$. If $\alpha < -1/q$, then
$$(L_1, L_2)_{1,q,\alpha} \cong L^{\mathcal{R}}_{2,q,\alpha,2}, \qquad (l_\infty, l_2)_{1,q,\alpha} \cong l^{\mathcal{L}}_{2,q,\alpha,2}.$$

If $\alpha > -1/q$, then
$$(L_2, L_\infty)_{0,q,\alpha} \cong L^{\mathcal{L}}_{2,q,\alpha,2}, \qquad (l_2, l_1)_{0,q,\alpha} \cong l^{\mathcal{R}}_{2,q,\alpha,2}.$$

By Lemmas 2.5 and 2.8 we get following embeddings (sf. [14], Corollaries 7.4 and 7.8).

**Lemma 2.9.** Let $0<q\leq\infty$. If $\alpha<-1/q$ (or $\alpha>-1/q$), then
$$L_{2,q}(\log L)_{\alpha+1/\min(q,2)} \subset L^{\mathcal{R}}_{2,q,\alpha,2} \text{ (or } L^{\mathcal{L}}_{2,q,\alpha,2}, \text{ resp.)} \subset$$
$$\subset L_{2,q}(\log L)_{\alpha+1/\max(q,2)} \cap L_{2,\max(q,2)}(\log L)_{\alpha+1/q}, \tag{2.6}$$

$$l_{2,q}(\log l)_{\alpha+1/\min(q,2)} \subset l^{\mathcal{L}}_{2,q,\alpha,2} \text{ (or } l^{\mathcal{R}}_{2,q,\alpha,2}, \text{ resp.)} \subset$$
$$\subset l_{2,q}(\log l)_{\alpha+1/\max(q,2)} \cap l_{2,\max(q,2)}(\log l)_{\alpha+1/q}. \tag{2.7}$$

**Remark 2.10.** For $0<q<2$ the spaces $L_{2,q}(\log L)_{\alpha+1/\max(q,2)}$ and $L_{2,\max(q,2)}(\log L)_{\alpha+1/q}$ are incomparable [37].

**Corollary 2.11.** If $\alpha<-1/2$ (or $\alpha>-1/2$) we have isomorphisms:
$$L^{\mathcal{R}}_{2,2,\alpha,2} \text{ (or } L^{\mathcal{L}}_{2,2,\alpha,2}, \text{ resp.)} \cong L_{2,2}(\log L)_{\alpha+1/2}, \quad l^{\mathcal{L}}_{2,2,\alpha,2} \text{ (or } l^{\mathcal{R}}_{2,2,\alpha,2}, \text{ resp.)} \cong l_{2,2}(\log l)_{\alpha+1/2}. \tag{2.8}$$

For the scale $l^{\mathcal{L}}_{2,q,\alpha,2}$ we also need the following inclusion.

**Lemma 2.12.** If $0<q<\infty$ and $\alpha<-1/q$, then $l^{\mathcal{L}}_{2,q,\alpha,2} \subset l^{\mathcal{L}}_{2,\infty,\alpha+\frac{1}{q},2}$.

**Proof.** Let $\{c_k\} \in l^{\mathcal{L}}_{2,q,\alpha,2}$. This means that

$$\| \{c_k\} \mid l^{\mathcal{L}}_{2,q,\alpha,2}\| = \left(\sum_{k=1}^{\infty}\left[(1+\log k)^{\alpha}\left(\sum_{i=1}^{k}\left[c_i^*\right]^2\right)^{1/2}\right]^q k^{-1}\right)^{1/q} < \infty.$$

The sequence $\{b_k\} := \left(\sum_{i=1}^{k}\left[c_i^*\right]^2\right)^{q/2}$ is non-negative and non-decreasing. Because $q<\infty$ and $q\alpha<-1$, for $m \geq 1$ we obtain

$$(1+\log m)^{1+q\alpha} \approx \int_{m}^{\infty}(1+\log x)^{q\alpha}\frac{dx}{x} \approx \sum_{k=m}^{\infty}(1+\log k)^{q\alpha}\frac{1}{k}.$$

Therefore

$$(1+\log m)^{1+q\alpha} b_m \approx b_m \sum_{k=m}^{\infty}(1+\log k)^{q\alpha}\frac{1}{k} \leq \sum_{k=m}^{\infty}(1+\log k)^{q\alpha} b_k \frac{1}{k}.$$

So, for all $m \geq 1$

$$(1+\log m)^{\alpha+\frac{1}{q}}\left(\sum_{i=1}^{m}\left[c_i^*\right]^2\right)^{1/2} \prec \left(\sum_{k=m}^{\infty}\left[(1+\log k)^{\alpha}\left(\sum_{i=1}^{k}\left[c_i^*\right]^2\right)^{1/2}\right]^q k^{-1}\right)^{1/q},$$

and hence

$$\|\{c_k\}\mid l^{\mathcal{L}}_{2,\infty,\alpha+\frac{1}{q},2}\| = \sup_{k\geq 1}(1+\log k)^{\alpha+\frac{1}{q}}\left(\sum_{i=1}^{k}\left[c_i^*\right]^2\right)^{1/2} \prec$$

$$\prec \left(\sum_{k=1}^{\infty}\left[(1+\log k)^{\alpha}\left(\sum_{i=1}^{k}\left[c_i^*\right]^2\right)^{1/2}\right]^q k^{-1}\right)^{1/q} = \|\{c_k\}\mid l^{\mathcal{L}}_{2,q,\alpha,2}\|. \quad \square$$

## 3. Proof of Theorem 1.1, corollaries, and remarks

**Proof of Theorem 1.1.** The assertion of Theorem 1.1 can be reformulated as follows: if $0<q\leq\infty$ and $\alpha< -1/q$ than the Fourier series map $\mathcal{F}$ is bounded from $L^{\mathcal{R}}_{2,q,\alpha,2}$ to $l^{\mathcal{L}}_{2,q,\alpha,2}$ and

$$\| \mathcal{F} \mid L^{\mathcal{R}}_{2,q,\alpha,2} \to l^{\mathcal{L}}_{2,q,\alpha,2} \| \prec \min(M, (1+\log M)^{|\alpha|}). \tag{3.1}$$

By Lemma 2.8 we have $(L_1, L_2)_{1,q,\alpha} \cong L^{\mathcal{R}}_{2,q,\alpha,2}$ and $(l_\infty, l_2)_{1,q,\alpha} \cong l^{\mathcal{L}}_{2,q,\alpha,2}$. Now due to (1.2) and (1.3) it only remains to apply Lemma 2.2. $\quad\square$

**Remark 3.1.** For the system $\{e^{i2\pi nx}\}$ Theorem 1.1 recovers [14, Theorem 8.2 (b)].

Using isomorphisms (2.8) we get the following corollary.

**Corollary 3.2.** If $\alpha < -1/2$, then

$$\|\mathcal{F} \mid L_{2,2}(\log L)_{\alpha+1/2} \to l_{2,2}(\log l)_{\alpha+1/2}\| \prec \min(M, (1+\log M)^{|\alpha|}).$$

**Remark 3.3.** The special case $\alpha = -1$ and the system $\{e^{i2\pi nx}\}$ of Corollary 3.2 recovers [11, Theorem 8.5].

The following result is a consequence of (3.1) and the left side of the inclusions (2.6).

**Corollary 3.4.** If $2 \leq q \leq \infty$ and $\alpha < 1/2 - 1/q$, then

$$\|\mathcal{F} \mid L_{2,q}(\log L)_\alpha \to l^{\mathcal{L}}_{2,q,\alpha-1/2,2}\| \prec \min(M, (1+\log M)^{|\alpha-1/2|}).$$

In particular, due to Remark 2.3, if $2 < q \leq \infty$, then

$$\|\mathcal{F} \mid L_{2,q} \to l^{\mathcal{L}}_{2,q,-1/2,2}\| \prec (1+\log M)^{1/2}.$$

**Remark 3.5.** In [5 and 6] (see also [29]) Bochkarev has proven the following estimate

$$\|\mathcal{F} \mid L_{2,q} \to l^{\mathcal{L}}_{2,\infty,1/q-1/2,2}\| \prec M \qquad (2 < q \leq \infty). \qquad (3.2)$$

Note that in the case $q=\infty$ formula (3.2) was proven in [32]. In [30] Bochkarev's inequality was improved in Lorentz–Zygmund spaces:

$$\|\mathcal{F} \mid L_{2,q}(\log L)_\alpha \to l^{\mathcal{L}}_{2,\infty,1/q-1/2+\alpha,2}\| \prec M \qquad (2<q\leq\infty, \alpha<1/2-1/q). \quad (3.3)$$

Due to Lemma 2.12 if $2<q<\infty$ and $\alpha <1/2 -1/q$ we have $l^{\mathcal{L}}_{2,q,\alpha-1/2,2} \subset l^{\mathcal{L}}_{2,\infty,1/q-1/2+\alpha,2}$. So, Corollary 3.4 improves both inequalities (3.2) and (3.3).

We finish this section by applying Lemma 2.9 to both spaces in (3.1).

**Corollary 3.6.** Let $0<q\leq\infty$ and $\alpha<-1/q$; then

$$\|\mathcal{F}\mid L_{2,q}(\log L)_{\alpha+1/\min(q,2)} \to l_{2,q}(\log l)_{\alpha+1/\max(q,2)} \cap l_{2,\max(q,2)}(\log l)_{\alpha+1/q}\| \prec \min(M,(1+\log M)^{|\alpha|}).$$

**Remark 3.7.** For the system $\{e^{i2\pi nx}\}$ Corollary 3.6 recovers [7, Theorem 5.3].

## 4. Proof of Theorem 1.2, corollaries and remarks

**Proof of Theorem 1.2.** The assertion of Theorem 1.2 can be reformulated as follows: if $0<q\leq\infty$ and $\alpha> -1/q$, then the operator $\mathcal{F}^{-1}$ is bounded from $l^{\mathcal{R}}_{2,q,\alpha,2}$ to $L^{\mathcal{L}}_{2,q,\alpha,2}$ and for any $\varepsilon>0$

$$\|\mathcal{F}^{-1} \mid l^{\mathcal{R}}_{2,q,\alpha,2} \to L^{\mathcal{L}}_{2,q,\alpha,2}\| \prec \min(M, (1+\log M)^{|\alpha|+\varepsilon+1/q}). \qquad (4.1)$$

By Lemma 2.8 we have $(l_2, l_1)_{0,q,\alpha} \cong l^{\mathcal{R}}_{2,q,\alpha,2}$ and $(L_2, L_\infty)_{0,q,\alpha} \cong L^{\mathcal{L}}_{2,q,\alpha,2}$. Now due to (1.4) and (1.5) it only remains to apply Lemma 2.2. $\square$

The following result is a consequence of (4.1) and the left side of the inclusions (2.7).

**Corollary 4.1.** Let $0<q\leq\infty$ and $\alpha>-1/q$; then for any $\varepsilon>0$

$$\|\mathcal{F}^{-1} \mid l_{2,q}(\log l)_{\alpha+1/\min(q,2)} \to L^{\ell}_{2,q,\alpha,2}\| \prec \min(M, (1+\log M)^{|\alpha|+\varepsilon+1/q}).$$

Applying Lemma 2.9 to the both spaces in (4.1) we get the next corollary.

**Corollary 4.2.** Let $0<q\leq\infty$ and $\alpha>-1/q$; then for any $\varepsilon>0$

$$\|\mathcal{F} \mid l_{2,q}(\log l)_{\alpha+1/\min(q,2)} \to L_{2,q}(\log L)_{\alpha+1/\max(q,2)} \cap L_{2,\max(q,2)}(\log L)_{\alpha+1/q}\| \prec$$

$$\prec \min(M,(1+\log M)^{|\alpha|+\varepsilon+1/q}).$$

In particular, if $\alpha>-1/2$, then for any $\varepsilon>0$

$$\|\mathcal{F}^{-1} \mid l_{2,2}(\log l)_{\alpha+1/2} \to L_{2,2}(\log L)_{\alpha+1/2}\| \prec \min(M, (1+\log M)^{|\alpha|+\varepsilon+1/2}).$$

**Remark 4.3.** For the system $\{e^{i2\pi nx}\}$ Theorem 1.2 and Corollary 4.2 recovers [14, Theorem 8.2 (g) and Lemma 8.4] respectively.

Using Remark 2.3 and Corollaries 4.1 and 4.2 we can formulate following assertions.

**Corollary 4.4.** Let $1<q\leq\infty$ and $\max(-1/q, 1/q-1)<\alpha<1-1/q$; then for any sufficiently small $\varepsilon>0$

$$\|\mathcal{F}^{-1} \mid l_{2,q}(\log l)_{\alpha+1/\min(q,2)} \to L^{\ell}_{2,q,\alpha,2}\| \prec (1+\log M)^{|\alpha|+\varepsilon+1/q}$$

and

$$\|\mathcal{F}^{-1} \mid l_{2,q}(\log l)_{\alpha+1/\min(q,2)} \to L_{2,q}(\log L)_{\alpha+1/\max(q,2)} \cap L_{2,\max(q,2)}(\log L)_{\alpha+1/q}\| \prec (1+\log M)^{|\alpha|+\varepsilon+1/q}.$$